\documentclass[journal]{IEEEtran}

\usepackage{algorithmic}
\usepackage{algorithm}
\usepackage{amsmath}
\usepackage{amssymb}
\usepackage{amsthm}
\usepackage{graphicx}
\usepackage{hyperref}
\usepackage{multirow}

\DeclareMathOperator*{\diag}{diag}

\DeclareMathOperator*{\minimize}{minimize}

\DeclareMathOperator*{\tr}{tr}
\DeclareMathOperator*{\vect}{vec}

\begin{document}

\title{A Fast Algorithm for Sparse Controller Design}
\author{Matt~Wytock,~J.~Zico~Kolter
\thanks{M. Wytock and J.Z. Kolter are with the School of Computer
  Science, Carnegie Mellon University, Pittsburgh, PA 15213 USA
  e-mail:\url{mwytock@cs.cmu.edu}, \url{zkolter@cs.cmu.edu}}}

\maketitle

\begin{abstract}
We consider the task of designing sparse control laws for large-scale
systems by directly minimizing an infinite horizon quadratic cost with
an $\ell_1$ penalty on the feedback controller gains. Our focus is on
an improved algorithm that allows us to scale to large
systems (i.e. those where sparsity is most useful) with convergence
times that are several orders of magnitude faster than existing
algorithms. In particular, we develop an efficient proximal Newton
method which minimizes per-iteration cost with a coordinate descent
active set approach and fast numerical solutions to the Lyapunov
equations. Experimentally we demonstrate the appeal of this approach
on synthetic examples and real power networks significantly larger
than those previously considered in the literature.
\end{abstract}

\section{Introduction}

This paper considers the task of designing sparse (i.e.,
decentralized) linear control laws for large-scale linear systems.
Sparsity and decentralized control have a long history
in the control literature: unlike the centralized control setting,
which in the $\mathcal{H}_2$ and $\mathcal{H}_\infty$ settings can be
solved optimally \cite{anderson1990optimal}, it has been known for
some time that the task of finding an optimal control law with
structure constraints in a hard problem \cite{blondel2000survey}.
Witenhausen's counterexample \cite{witsenhausen1968counterexample}
famously demonstrated that even for a simple linear system, a linear
control law is no longer optimal, and subsequent early work focused
on finding effective decentralized controllers for specific problem
instances \cite{lau1972decentralized} or determining the theoretical
existence of stabilizing distributed control methods
\cite{wang1973stabilization}; they survey of Sandell et al.,
\cite{sandell1978survey} covers many of these earlier approaches in
detail.

In recent years, there has been an increasing interest in sparse and
decentralized control methods, spurred 1) by increasing interest in
large-scale systems such as the electrical grid, where some form of
decentralization seems critical for practical control strategies, and 2)
by increasing computational power that can allow for effective
controller design methods in such systems.  Generally, this work has
taken one of two directions.  On the one hand, several authors have
looked at restricted classes of dynamical systems where the
true optimal control law is provably sparse, resulting in
efficient methods for computing optimal decentralized controllers
\cite{fan1994centralized,qi2004structured,rotkowitz2006characterization,shah2010H2}.
A notable recent example of such work has been the
characterization of all systems that admit convex constraints on the
controller (and thus allow for exact sparsity-constrained controller
design) using the notion of quadratic invariance
\cite{rotkowitz2006characterization,swigart2011decentralized}.  On the
other hand, an alternative approach has been to search for
\emph{approximate} (suboptimal) decentralized controllers, either by
directly solving a nonconvex optimization problem
\cite{lin2013design}, by constraining the class of allowable Lyapunov
functions in a convex parameterization of the optimal $\mathcal{H}_2$
or $\mathcal{H}_\infty$ controllers
\cite{schuler2011design,schuler2013decentralized} or by employing a
convex, alternative optimization objective as opposed to the typical
infinite horizon cost
\cite{dvijotham2013convexity,dvijotham2013convex}.  This present paper
follows upon this second line of work, specifically building upon the
framework established in \cite{lin2013design}, which uses $\ell_1$
regularization (amongst other possible regularizers) to
\emph{discover} the good sparsity patterns in the control law (though
our method also applies directly to the case of a fixed sparsity
pattern).

Despite the aforementioned work, an element that has been notably
missing from past work in the area is a focus on the
\emph{algorithmic} approaches that can render these methods practical
for large-scale systems, such as those with thousands of states and
controls or more.  Indeed, it is precisely for such systems that
sparse and decentralized control is most appealing, and yet most past
work we are aware of has focused solely on semidefinite programming
formulation of the resulting optimization problems
\cite{schuler2011design,schuler2013decentralized} (which scale poorly
in off-the-shelf solvers), or very approximate first order or
alternating minimization methods \cite{lin2013design} such as the
alternating direction method of multipliers
\cite{boyd2011distributed}.  As a result, most of the demonstrated
performance of the methods in these past papers has focused solely on
relatively small-scale systems.  In contrast, sparse methods in
fields like machine learning and statistics, which have received a
great deal of attention in recent years
\cite{tibshirani1996regression,donoho2006compressed,candes2006stable},
evolved simultaneously with efficient algorithms for solving these
statistical estimation problems
\cite{efron2004least,koh2007interior}.  The goal of the present paper
is to push this algorithmic direction in the area of sparse control,
developing methods that can make handle large-scale sparse controller
design.

\subsection{Contributions of this paper}

In this paper, we develop a fast algorithm for large-scale sparse
state-feedback controller design that is several orders of magnitude
faster than previous approaches.  The approach is based upon a second
order method known as a Newton-Lasso approach (also called proximal
Newton methods) \cite{tseng2009coordinate,byrd2013inexact},
which iteratively performs \emph{regularized} Newton
steps to minimize a smooth objective plus an $\ell_1$ penalty.

As with any Newton method, although the number of iterations can be
much much less than for first-order methods, the danger is that the
\emph{per-iteration} cost can be substantially higher, to the point
that the overall algorithm is in fact slower for larger problems.
This is of particular concern in the optimal control setting, where,
as we show below, computing a \emph{single} inner product with the
Hessian typically involves solving a Lyapunov equation in $n^2$
variables, which itself is an $O(n^3)$ operation.  Thus, the majority
of the algorithmic work involves developing a method where each Newton
step has relatively low cost, so that the overall speed of the
algorithm is indeed much faster than alternative approaches.

We accomplish these speedups in two ways. First, we employ a
coordinate descent active set approach to solve each (regularized)
Newton step in the algorithm.  These methods are particularly
appealing for $\ell_1$ regularized problems, as coordinate descent
approaches typically work quite well for many quadratic $\ell_1$
problems \cite{tseng2009coordinate}, and thus are well-suited to
solving the inner Newton steps.  In past work, such methods have been
applied with great success to $\ell_1$ regularized problems
such as sparse inverse covariance estimation
\cite{hsieh2013sparse} and spare conditional Gaussian
models \cite{wytock2013sparse}; indeed, these methods currently
comprise the state of the art for solving these optimization
problems.  Second, although the general Newton coordinate descent
framework looks appealing for such $\ell_1$ regularized problems, its
application to the sparse controller design setting is challenging: a
naive implementation of coordinate descent requires many Hessian inner
products per outer iteration, and as we saw above, each such inner
product is very costly for this setting.  Thus, a significant
component of the algorithm relates to how we can reduce the
cost of these inner products from $O(n^3)$ to $O(n)$ through a series
of precomputations and transformations.  To do this, we make use of
methods from the numerical solution of Lyapunov equations
\cite{simoncini2013computational}, the fast multipole method
\cite{greengard1987fast}, and the Autonne-Kakagi factorization
\cite[Corollary 2.6.6]{horn2012matrix}.

The end result is a method that achieves high numerical precision while
being in many cases orders of magnitude faster than existing
approaches.  This represents a substantial advance in the practical
application of such approaches, allowing them to be applied to
problems that previously had been virtually unsolvable in the sparse
controller design framework.

\section{The sparse optimal control framework}

Here we formally define our control and optimization framework, based
upon the setting in \cite{lin2013design}.  Formally, we consider the
linear Gaussian system
\begin{equation}
\dot{x}(t) = A x(t) + B u(t) + W^{1/2}\epsilon(t)
\end{equation}
where $x \in \mathbb{R}^n$ denotes the state variables, $u \in
\mathbb{R}^m$ denotes the control inputs, $\epsilon$ is a zero-mean
Brownian motion process, $A \in \mathbb{R}^{n \times
  n}$ and $B \in \mathbb{R}^{n \times m}$ are system matrices, and
$W^{1/2} \in \mathbb{R}^{n \times n}$ is a noise covariance matrix.
We  seek to optimize the infinite horizon LQR cost for a linear
state-feedback control law $u(t) = Kx(t)$ for $K \in \mathbb{R}^{m
  \times n}$,
\begin{equation}
\label{eq-cost}
J(K) = \lim_{T \rightarrow \infty} \frac{1}{T} \mathbf{E} \left [\int_0^T
  \left( x(t)^T Q x(t) + u(t)^T R u(t) \right) dt \right ]
\end{equation}
which can also be written in the alternative form
\begin{equation}
J(K) = \left \{ \begin{array}{ll} \tr PW = \tr L(Q + K^TRK) &
  \mbox{ $A+BK$ stable}\\
\infty &\mbox{ otherwise}. \end{array} \right .
\end{equation}
where $L = L(K)$ and $P = P(K)$ are the unique solutions to the Lyapunov
equations
\begin{equation}
\begin{split}
& (A + BK)L + L(A + BK)^T + W = 0 \\
& (A + BK)^TP + P(A+BK) + Q + K^T R K = 0.
\end{split}
\end{equation}
When $K$ is unconstrained, it is well-known that the problem can
be solved by the classical LQR algorithm, though this results in a
dense (i.e., centralized) control law, where each control will
typically depend on each state.

To encourage sparsity in the controller, \cite{lin2013design} proposed
to add an additional penalty to the $\ell_1$ norm of the controller
$K$ (along with other possible regularization terms that we don't consider
directly here).  Here we will consider this framework with a weighted
$\ell_1$ norm: i.e., we are concerned with solving the optimization
problem
\begin{equation}
\minimize_K J(K) + g(K)
\end{equation}
where $J(K)$ defined in \eqref{eq-cost} is the LQR cost (or the
$\mathcal{H}_2$ norm for output $z(t) = (Q^{1/2} x(t), R^{1/2} u(t))$ and
treating $\epsilon(t)$ as a disturbance input) and $g(K)$ is the
sparsity-promoting penalty which we take to be
\begin{equation}
g(K) = \|\Lambda \circ  K \|_1 = \sum_{ij} \Lambda_{ij}|K_{ij}|,
\end{equation}
a weighted version of the $\ell_1$ norm.  This formulation allows us
to use this single algorithmic framework to capture both traditional
$\ell_1$ regularization of the control matrix, as well as optimization
using a fixed pattern of nonzeros, by setting the approximate
elements of $\Lambda$ to $0$ or $\infty$.

Our algorithm uses gradient and Hessian information extensively.
Following standard results \cite{rautert1997computational}, the
gradient of $J(K)$ is given by
\begin{equation}
\nabla J(K) = 2(B^TP + RK)L
\end{equation}
when $A + BK$ is stable.  The Hessian is somewhat cumbersome to
formulate directly, but we can concisely write its inner
product with a direction $D \in \mathbb{R}^{m \times n}$ as
\begin{equation}
\label{eq-hessian-inner}
\begin{split}
& \vect(D)^T\nabla^2J(K)\vect(D) = \\
& \quad\quad 2 \tr \left( \tilde{L}(PB + K^TR) + L(\tilde{P}B + D^T)R
\right) D
\end{split}
\end{equation}
where $\vect(\cdot)$ denotes the vectorization of a matrix, and
$\tilde{L} = \tilde{L}(K,D)$ and $\tilde{P} = \tilde{P}(K,D)$
are the unique solutions to two more Lyapunov equations
\begin{equation}
\label{eq-lyap-diff}
\begin{split}
& (A + BK)\tilde{L}  + \tilde{L}(A + BK)^T + BD L  + LD^TB^T = 0 \\
& (A + BK)^T\tilde{P} + \tilde{P}(A + BK) + ED + D^TE^T = 0,
\end{split}
\end{equation}
denoting $E = PB + K^TR$ for brevity.  Since solving these Lyapunov
equations takes $O(n^3)$ time, evaluating the function and gradient,
or evaluating a single inner product with the Hessian, are all
$O(n^3)$ operations.

Traditionally, direct second order Newton methods have seen relatively
little application in Lyapunov-based control, precisely because
computing these Hessian terms is computationally intensive.  Instead,
typical approaches have focused on approaches that use gradient
information only, either in a quasi-Newton setup
\cite{rautert1997computational}, or by including only certain terms
from the Hessian as in the Anderson-Moore method
\cite{anderson1990optimal}.  However, since a single iteration of any
first order approach is already reasonably expensive $O(n^3)$
operation, the significantly reduced iteration count of typical Newton
methods looks appealing, provided we have a way to efficiently compute
the Newton step.  The algorithm we propose in the next section does
precisely that, bringing down the complexity of a Newton step to a
computational cost similar to that of just evaluating the function.

\section{A fast Newton-Lasso algorithm}
\label{sec-algorithm}

\subsection{Overview of the algorithm}

\begin{algorithm}[t]
\caption{Newton-CD for Sparse LQR}
\label{alg-newton-lasso}
\begin{algorithmic}
\STATE \textbf{Input:} Stochastic linear system $A$, $B$, $W$; Regularization
parameters $Q$, $R$, $\Lambda$
\STATE \textbf{Output:} Optimal controller $K$
\STATE \textbf{Initialize:} $K \leftarrow 0$
\WHILE{(not converged)}
  \STATE 1. Compute the active set $\mathcal{A}$
\[
\mathcal{A} = \{(i,j) \mid |(\nabla J(K))_{ij}| > \Lambda_{ij} \text{ or } K_{ij} \ne 0 \}
\]
  \STATE 2. Find the regularized Newton step using coordinate descent over the
  active set $\mathcal{A}$
  \[
\hat{D} = \arg \min_{D \in \mathcal{D_\mathcal{A}}}  \tilde{J}_K(D) + \|(K + D) \circ \Lambda\|_1
\]
  \STATE 3. Choose a stepsize $\alpha$ and set $K \leftarrow K + \alpha \hat{D}$
\ENDWHILE
\end{algorithmic}
\end{algorithm}

Our algorithm follows the overall structure of a ``Newton-Lasso''
method \cite{tseng2009coordinate}, also sometimes called proximal
Newton methods \cite{byrd2013inexact}. The overall idea is to
repeatedly form a second order approximation to the smooth  component
of the objective function, $J(K)$, and minimize
this second order approximation plus an $\ell_1$ regularization term.
This effectively reduces the problem from solving an arbitrary
smooth objective with $\ell_1$ regularization to a \emph{quadratic}
objective with $\ell_1$ regularization, which is typically an
easier problem in practice (though here it still requires iterative
optimization itself).  But, as with Newton's method, the number of
``outer loop'' iterations needed is typically very small.

Of course, the overall time complexity of the algorithm depends
critically on the efficiency of the inner loop for finding the
regularized Newton direction, a step which cannot be computed in
closed form.  In this work, we specifically propose to use a coordinate
descent approach for finding these approximate Newton directions,
leading to an approach we refer to as Newton Coordinate Descent, or
simply Newton-CD.  Coordinate descent, though a simple method, is
known to perform quite well for quadratic $\ell_1$ regularized
problems \cite{tseng2009coordinate,hsieh2013sparse,wytock2013sparse}, but it's
real
benefit in our setting comes from two characteristics:
\begin{enumerate}
\item Coordinate descent methods allow us to optimize only over a
  relatively small ``active set'' $\mathcal{A}$ of size $k \ll mn$,
  which includes only the nonzero elements of $K$ plus elements
  with large gradient values.  For problems that exhibit substantial
  sparsity in the solution, this often lets us optimize over much
  fewer elements than would be required if we considered all the
  elements of $K$ at each iteration.
\item Most importantly, by properly precomputing certain terms,
  caching intermediate products, and exploiting problem structure, we
  can reduce the per-coordinate-update computation in coordinate
  descent from $O(n^3)$ (the naive solution, since each coordinate
  update requires computing an inner product with the Hessian matrix)
  to $O(n)$.
\end{enumerate}
We describe each of these two elements in detail in the following two
sections, beginning with a general description of the Newton-CD
approach using an active set and then detailing the computations
needed for efficient coordinate descent updates. A C++ and MATLAB implementation
is available at \url{http://www.cs.cmu.edu/~mwytock/lqr/}.




\subsection{The active set Newton-CD algorithm}
The generic Newton-CD algorithm is shown in Algorithm
\ref{alg-newton-lasso}.   Formally, at each iteration we form the
second order Taylor expansion
\begin{equation}
\begin{split}
J(K + D)
&\approx \tr \nabla
J(K)^TD + \frac{1}{2}\vect(D)^T\nabla^2J(K)\vect(D) \\
&\equiv \tilde{J}_K(D)
\end{split}
\end{equation}
which in our problem takes the form
\begin{equation}
\tilde{J}_K(D) = 2\tr LED + \tr \tilde{L}ED + \tr L\tilde{P}BD +
\tr LD^TRD
\end{equation}
where $\tilde{L}$ and $\tilde{P}$ are defined implicitly as the unique
solutions to the Lyapunov equations given in \eqref{eq-lyap-diff}. We
minimize this quadratic function over the active set $\mathcal{A}$,
only including coordinates that are nonzero at the current iterate or
are nonoptimal---those that satisfy
\begin{equation}
|(\nabla J(K))_{ij}| > \Lambda_{ij} \text{ or } K_{ij} \ne 0
\end{equation}
defining $\mathcal{D_\mathcal{A}}
\equiv \{ D : D_{ij} = 0 \forall (i,j) \notin \mathcal{A} \}$, the set of
candidate directions. Intuitively, these correspond to only those
elements that are either non-zero or violate the current optimality
conditions of the optimization problem.

In contrast to the standard Newton method, the Newton-Lasso
then minimizes the second order approximation with the addition of
(weighted) $\ell_1$ regularization
\begin{equation}
\label{eq-argmin-dir}
\hat{D} = \arg \min_{D \in \mathcal{D_\mathcal{A}}} \tilde{J}_K(D) +
\|(K + D) \circ \Lambda\|_1
\end{equation}
and updates $K \leftarrow K + \alpha \hat{D}$ where
$\alpha$ is chosen using backtracking line search and an Armijo
rule. Intuitively, it can be shown that $\alpha \to 1$ as the
algorithm progresses, causing the weighted $\ell_1$ penalty on $K + D$
to shrink $D$ in the direction the promotes sparsity in $K + \hat{D}$;
the weights $\Lambda$ control how aggressively we shrink each
coordinate: $\Lambda_{ij} \to \infty$ forces $(K + \hat{D})_{ij} \to
0$.  Importantly, we note that with this formulation we will
never choose a $K$ such that $A+BK$ is unstable; this would make the
resulting objective infinite, and a smaller step size would be
preferred by the backtracking line search.

As mentioned, in order to find the regularized Newton step
efficiently, we use coordinate descent which is appealing for Lasso
problems as each coordinate update can be computed in closed
form. This reduces \eqref{eq-argmin-dir} to iteratively
minimizing each coordinate
\begin{equation}
\label{eq-cd}
\hat{\mu} = \arg\min_\mu \tilde{J}_K(D + \mu e_ie_j^T) + \Lambda_{ij}|K_{ij} + D_{ij} + \mu|
\end{equation}
where $e_i$ denotes the $i$th basis vector, and then setting $D
\leftarrow D + \hat{\mu} e_ie_j^T$. Since the second
order approximation $\tilde{J}_K(D)$ is quite complex (it
depends on the solution to four Lyapunov equations, two of which
depend on $D$), deriving efficient coordinatewise updates is somewhat
involved. In the next section we describe how each coordinate descent
iteration can be computed in $O(n)$ time by precomputing a single
eigendecomposition and solving the Lyapunov equations explicitly.

\subsection{Fast coordinate updates}
\begin{algorithm}[t]
\caption{Coordinate descent for regularized Newton step}
\label{alg-cd}
\begin{algorithmic}
\STATE \textbf{Input:} Stochastic linear system $A$, $B$, $W$; regularization
parameters $Q$, $R$, $\Lambda$; current iterate $K$; active set $\mathcal{A}$;
solution to Lyapunov equations $L, P$
\STATE \textbf{Output:} Regularized Newton step $D$
\STATE \textbf{Initialize}: $D \leftarrow 0$, $\Psi \leftarrow 0$
\STATE 1. Compute the eigendecomposition $A + BK = USU^{-1}$
\STATE 2. Let $\Theta_{ij} = 1/(S_{ii} + S_{jj})$ and compute $\Theta = XX^T$
\STATE 3. Precompute matrix products as in \eqref{eq-coord-precompute}
\WHILE{(not converged)}
  \FOR{coordinate $(i,j)$ in  $\mathcal{A}$}
  \STATE 1. Compute $a,b,c$ according to \eqref{eq-coord-updates}
  and set
  \begin{equation*}
    \mu = -c + S_{\lambda/a}\left( c - \frac{b}{a} \right)
  \end{equation*}
  \STATE 2. Update solution
  \begin{equation*}
    D_{ij} \leftarrow D_{ij} + \mu
  \end{equation*}
  \STATE 3. Update the cached matrix products
  \begin{equation*}
    \begin{split}
      (\Psi^0)_j &\leftarrow (\Psi^0)_j + \mu R_i \\
      (\Psi_k^1)_i &\leftarrow (\Psi_k^1)_i + \mu (\Phi_k^1)_j \\
      (\Psi_k^2)_j &\leftarrow (\Psi_k^2)_j + \mu (\Phi_k^4)_i \\
      (\Psi_k^3)_j &\leftarrow (\Psi_k^3)_j + \mu (\Phi_k^2)_i \\
      (\Psi_k^4)_j &\leftarrow (\Psi_k^4)_j + \mu (\Phi_k^3)_i
    \end{split}
  \end{equation*}
  where (in general) $A_j$ denotes the $j$th column of $A$
  \ENDFOR
\ENDWHILE
\end{algorithmic}
\end{algorithm}

To begin, we consider the explicit forms for $\tilde{L}(K,D)$ and
$\tilde{P}(K,D)$, the unique solutions to the Lyapunov equations
depending on $D$. Since each coordinate descent update changes an
element of $D$, a naive approach would require re-solving these two
Lyapunov equations and $O(n^3)$ operations per iteration. Instead,
assuming that $A + BK$ is diagonalizable, we precompute a single
eigendecomposition $A + BK = USU^{-1}$ and use  this to compute the
solutions to the Lyapunov equations directly. For example, the
equation describing $\tilde{L}$ can be written as
\begin{equation}
\begin{split}
USU^{-1}\tilde{L}  + \tilde{L}U^{-T}SU^T + BDL  + LD^TB^T = 0
\end{split}
\end{equation}
and pre and post multiplying by $U^{-1}$ and $U^{-T}$ respectively gives
\begin{equation}
S\tilde{L}_U  + \tilde{L}_US  = -U^{-1}(BDL  + LD^TB^T)U^{-T}
\end{equation}
where $\tilde{L}_U = U^{-1}\tilde{L}U^{-T}$. Since $S$ is diagonal, this
equation has the solution
\begin{equation}
(\tilde{L}_U)_{ij} = -\frac{(U^{-1}(BDL  + LD^TB^T)U^{-T})_{ij}}{S_{ii} + S_{jj}}
\end{equation}
which we rewrite as the Hadamard product
\begin{equation}
\tilde{L}_U = U^{-1} (BDL + LD^TB^T) U^{-T} \circ \Theta
\end{equation}
with $\Theta_{ij} = -1/(S_{ii} + S_{jj})$.

\textbf{Fast $\Theta$ multiplication via the Fast Multiple Method}:
Precomputing the eigendecomposition of $A + BK$ in this manner
immediately allows for an $O(n^2)$ algorithm for evaluating Hessian
products, but reducing this to $O(n)$ requires exploiting additional
structure in the problem.  In particular, we consider the form of the
$\Theta$ matrix above, which is an example of a Cauchy matrix,
that is, matrices with the form $C_{ij} = 1/(a_i - b_j)$.
Like several other special classes of matrices, matrix-vector products
with a Cauchy matrix can be computed more quickly than for a standard
matrix.  In particular, the Fast Multiple Method (FMM)
\cite{greengard1987fast},
specifically the 2D FMM using the Laplace kernel, provides an $O(n)$
algorithm (technically $O(n \log\frac{1}{\epsilon})$ where $\epsilon$
is the desired accuracy) for computing the matrix vector product
between a Cauchy matrix and an arbitrary vector.\footnote{In theory,
  such matrix-vector products for Cauchy matrices can be computed
  exactly in time $O(n \log^2 n)$ \cite{gerasoulis1988fast}, but these
  approaches are   substantially less numerically robust than the FMM,
  so the FMM is typically preferred in practice \cite{pan1998new}.}

Although the FMM provides a theoretical method for quickly computing
Hessian inner products, in our setting the overhead
involved with actually setting up the factorization (which also takes
$O(n)$ time, but with a relatively larger constant) would make using
an off-the-shelf implementation of the FMM quite costly. However, our
setting in fact is somewhat easier as $\Theta$ is fixed per outer
Newton iteration; thus we can factor $\Theta$ once at (relatively)
high computation cost and then directly use this factorization is
subsequent iterations.  Each FMM operation implicitly factors $\Theta$
in a hierarchical manner with blocks of low-rank structure, though
here the situation is simpler: since we maintain $A+BK$ to
be stable at each iteration, all the eigenvalues are in
the left half plane and representing $\Theta$ as a Cauchy matrix
\begin{equation}
\Theta_{ij} = \frac{1}{a_i - b_j},
\end{equation}
i.e., $a_i = S_{ii}$, $b_j = -S_{jj}$ leads to points, $a_i,b_j \in
\mathbb{C}$ that are \emph{separated} in the context of the FMM.  This
means that $\Theta$ in fact simply admits a low-rank representation
(though the actual rank will be problem-specific, and depend on how
close the eigenvalues of $A+BK$ are to the imaginary axis).  Thus,
while slightly more advanced factorizations may be possible, for the
purposes of this paper we simply use the property, based upon the FMM,
that $\Theta$ will typically admit a low-rank factorization.

\textbf{The Autonne-Kagaki factorization of $\Theta$}:
Using the above property, we can compute the optimal low rank
factorization of $\Theta$ using the (complex) singular value
decomposition to obtain a factorization $\Theta = X Y^*$.  But since
$\Theta$ is a complex symmetric (but not Hermitian) matrix, it also can
be factored as $\Theta = V S V^T$ where $S$ is a diagonal matrix of
the singular values of $\Theta$ and $V$ is a complex unitary matrix
\cite[Corollary 2.6.6]{horn2012matrix}.
This factorization lets us speed up the resulting computations by 2
fold over simply using an SVD, as we have significantly fewer matrices
to precompute in the sequel.

Specifically, writing $\Theta = \sum_{i=1}^n x_i x_i^T$, and using the
fact that for a Hadamard product
\begin{equation}
A \circ ab^T = \diag(a)A\diag(b).
\end{equation}
we can write the Lyapunov solution $\tilde{L}$ analytically as
\begin{equation}
\tilde{L} = U \tilde{L}_U U^T = \sum_{k=1}^r X_i (BDL + (BDL)^T)X_i^T
\end{equation}
where we let $X_i = U \diag(x_i) U^{-1}$, the
transformed version of the diagonal matrix corresponding to the $i$th column of
$X$. With the same approach, we write the explicit form for $\tilde{P}$ as
\begin{equation}
\tilde{P} = \sum_{i=1}^r X_i^T((ED)^T + ED)X_i.
\end{equation}
Using these explicit forms for $\tilde{L}$ and $\tilde{P}$ we observe that
$\tr \tilde{L}ED = \tr \tilde{P}BDL$ and the second order Taylor expansion
simplifies to
\begin{equation}
\begin{split}
\tilde{J}_K(D) &= 2 \tr LED + \tr LD^TRD \\
&+ 2 \tr \sum_{i=1}^r X_i^T((ED)^T +
ED)X_iBDL.
\end{split}
\end{equation}

\textbf{Closed form coordinate updates}:
Next, we consider coordinatewise updates to minimize $\tilde{J}_K(D)$ with the
addition of $\ell_1$ regularization. In particular, consider optimizing over $\mu$
the rank one update $D + \mu e_ie_j^T$; for each term we get a quadratic
function in $\mu$ with coefficients that depend on several matrix products.
For example, the second term of $LD^TRD$ yields
\begin{equation}
\begin{split}
 & \tr L(D + \mu e_ie_j^T)^TR(D + \mu e_ie_j^T) \\
= &   \tr LD^TRD + 2 \mu(RDL)_{ij} + \mu^2R_{ii}L_{jj}.
\end{split}
\end{equation}
For each term in $\tilde{J}_K(D)$, we repeat these steps to derive
\begin{equation}
\begin{split}
 & \arg\min_\mu \tilde{J}_K(D + \mu e_ie_j^T) + \|(K + D + \mu e_ie_j^T) \circ \Lambda\|_1 \\
= & \arg\min_\mu \frac{1}{2} a\mu^2 + b \mu + \|c + \mu\|_1
\end{split}
\end{equation}
where
\begin{equation}
\begin{split}
a &= 2 R_{ii}L_{jj} \\
&+ 4 \left( \sum_{k=1}^r (E^TX_kB)_{ii}(LX_k^T)_{jj}
+ (LX^TE)_{ji}(XB)_{ji} \right) \\
b &= 2 (E^TL)_{ij} + 2(RDL)_{ij} \\
& + 2\left( \sum_{k=1}^r (E^TX_kBDLX_k^T)_{ij} +
(B^TX_k^TEDX_kL)_{ij} \right) \\
& + 2\left( \sum_{k=1}^r (X_kBDLX_k^TE)_{ji} + (LX_k^TEDX_kB)_{ji} \right) \\
c &= K_{ij} + D_{ij}.
\end{split}
\end{equation}
This has the closed form solution
\begin{equation}
  \mu = -c + S_{\lambda/a}\left( c - \frac{b}{a} \right)
\end{equation}
where $S_\lambda$ is the soft-thresholding operator.

\textbf{Caching matrices}:
Naive computation of these matrix products for $a$, $b$ and $c$ still requires
$O(n^3)$ operations; however, all matrices except $D$ remain fixed
over each iteration of the inner loop, allowing us to precompute many
matrix products. Let
\begin{equation}
\label{eq-coord-precompute}
\begin{split}
\Phi^0  &= LE \\
\Phi_k^1 &= X_kL \\
\Phi_k^2 &= X_kB \\
\Phi_k^3 &= LX_k^TE \\
\Phi_k^4 &= B^TX^TE.
\end{split}
\end{equation}
In addition as we iteratively update $D$, we also maintain the matrix
products
\begin{equation}
\begin{split}
\Psi^0   &= RD \\
\Psi_k^1 &= \Phi_k^1D^T \\
\Psi_k^2 &= \Phi_k^4D  \\
\Psi_k^3 &= \Phi_k^2D \\
\Psi_k^4 &= \Phi_k^3D
\end{split}
\end{equation}
which allows us to efficiently compute
\begin{equation}
\label{eq-coord-updates}
\begin{split}
a &= 2 R_{ii}L_{jj} + 4 \left( \sum_{k=1}^r (\Phi_k^4)_{ii}(\Phi_k^1)_{jj}
+ (\Phi_k^3)_{ij}(\Phi_k^2)_{ji} \right) \\
b &= 2 \left( (\Phi^0)_{ji} + (\Psi^0L)_{ij} \right) \\
& + 2\left( \sum_{k=1}^r (\Psi_k^1\Phi_k)_{ji} +
(\Psi_k^2\Phi_k^1)_{ij} \right) \\
& + 2\left( \sum_{k=1}^r (\Psi_k^3\Phi_k^3)_{ji} + \Psi_k^4\Phi_k^2)_{ji} \right) \\
c &= K_{ij} + D_{ij}
\end{split}
\end{equation}
resulting in an $O(n)$ time per iteration (in general we can compute an element
of a matrix product $(AB)_{ij}$ as the dot product between the $i$th row of $A$
and the $j$th column of $B$). Updating the cached products $\Psi$ also
requires $O(n)$ time as a change to a single coordinate of $D$ requires
modifying a single row or column of one of the products $\Psi$. The complete
algorithm is given in Algorithm \ref{alg-cd} and has $O(n)$  per iteration
complexity as opposed to the $O(n^3)$ naive implementation.

\subsection{Additional algorithmic elements}

While the above algorithm describes the basic second order approach,
several elements are important for making the algorithm practical and
robust to a variety of different systems.

\textbf{Initial conditions}:
One crucial element that affects the algorithm's performance is the
choice of initial $K$ matrix.  Since the objective $J(K)$ is infinite for
$A+BK$ unstable, we require that the initial value must
stabilize the system.  We could simply choose the full LQR
controller $K^{\mathrm{LQR}}$ as this initial point; it may take time
$O(n^3)$ to compute the LQR solution, but since our algorithm is
$O(n^3)$ overall, this is typically not a prohibitive cost.
However, the difficulty with this strategy is that the resulting
controller is not sparse, which leads to a full active set for the
first step of our Newton-CD approach, substantially slowing down the
method.  Instead, a single soft-thresholding step on the
LQR solution produces a good initial starting point that is both
guaranteed to be stable, and which leads to a much smaller active set
in practice. Formally, we compute
\begin{equation}
K^{(0)} = S_{\alpha\Lambda}(K^{\mathrm{LQR}})
\end{equation}
where $\alpha \le  1$ is chosen by backtracking line search such that
the regularized objective decreases and $K^{(0)}$ remains stable.

In addition, if the goal is to sweep across a large range of possible
regularization parameters, we can employ a ``warm start'' method that
initializes the controller to the solution of previous optimization
problems.

\textbf{Unstable initial controllers}:
In the event where we do not want to start at the LQR solution, it is
also possible to begin with some initial controller $K^{(0)}$ that is
\emph{not} stabilizing using a ``deflation'' technique.  Specifically,
rather than find an optimal control law for the linear system $(A,B)$,
we find a controller for the linear system $(A - \nu I,B)$ where $\nu$
is chosen such that $A - \nu I + BK^{(0)}$ is stable with some margin.
The resulting controller $K^{(1)}$ will typically stabilize the system
to a larger degree, and we can repeat this process until it
produces a stabilizing control law.  Further, we typically do not need
to run the Newton-CD method to convergence, but can often obtain a
better stabilizing control law after only a few outer iterations.

\textbf{Handling non-convexity}:
As mentioned above, the objective $J(K)$ is not a convex function, and
so can (and indeed, often does in practice) produce indefinite Hessian
matrices.  In such cases, the coordinate descent steps are not
guaranteed to produce a descent direction, and indeed can cause the
overall descent direction to diverge.  Furthermore, since we never
compute the full Hessian $\nabla_K^2 J(K)$ (even restricted to just
the active set), it is difficult to perform typical operations to
handle non-convexity such as projecting the Hessian onto the positive
definite cone.  Instead, we handle this non-convexity by a fallback to
a simpler quasi-Newton coordinate descent scheme
\cite{tseng2009coordinate}.  At each Newton
iteration we also form a coordinate descent update based upon a
\emph{diagonal} PSD approximation to the Hessian,
\begin{equation}
\bar{J}(K + D)=  \tr \nabla J(K)^T D +
\frac{1}{2}\vect(D)^T H \vect(D)
\end{equation}
where
\begin{equation}
H_{ii} = \min\{\max\{10^{-2}, (\nabla^2 J(K))_{ii}\}, 10^4\}.
\end{equation}
The diagonal terms of the Hessian are precisely the $a$ variables that
we compute in the coordinate descent iterations anyway, so this search
direction can be computed at little additional cost.  Then, we simply
perform a line search on both update directions simultaneously, and
choose the next iterate with the largest improvement to the
objective.  In practice, the algorithm sometimes uses the fallback
direction in the early iterations of the method until it converges to a
convex region around a (local) optimum where the full Newton step
causes much larger function decreases and so is nearly always
chosen.  This fallback procedure also provides a convergence
guarantee for our method: the quasi-Newton coordinate descent approach
was analyzed in \cite{tseng2009coordinate} and shown to converge for both
convex and non-convex objectives.  Since our algorithm always takes at
least as good a step as this quasi-Newton approach, the same
convergence guarantees hold here.

\textbf{Inner and outer loop convergence and approximation}:
Finally, a natural question that we don't address directly in the
above algorithmic presentations involves how many iterations (both
inner and outer), as well as what constitutes a sufficient
approximation for certain terms like the rank of $\Theta$.  In
practice, a strength of the Newton-CD methods is that it can be fairly
insensitive to slightly less accurate inner loops
\cite{byrd2013inexact}: in such cases, the approximate Newton
direction is still typically much better than a gradient direction,
and while additional outer loop iterations may be required, the
timing of the resulting method is somewhat insensitive to choice of
parameters for the inner loop convergence and for different
low-rank approximations to $\Theta$.  In our implementation, we run
the inner loop for at most $t/3$ iterations at outer loop iteration
$t$ or until the relative change in the direction $D$ is less than
$10^{-2}$ in the Frobenius norm.
\section{Experiments}
\label{sec-experiments}

In this section we evaluate the performance of the proposed algorithm on the
task of finding sparse optimal controllers for a synthetic mass-spring system
and wide-area control in power systems. In both
settings, the method finds sparse controllers that perform nearly
optimally while only depending on a small subset of the state
space; furthermore, as we scale to larger examples, we demonstrate that these
optimal controllers become \emph{more sparse}, highlighting the increased role
of sparsity in larger systems.

Computationally, we compare the
convergence rate of our algorithm
to that of existing approaches for solving the sparse optimal control problem
and demonstrate that the proposed method converges rapidly to highly
accurate solutions, significantly outperforming previous approaches. Although
the solution accuracy required for a ``good enough'' controller is
problem-specific, since iteration complexity grows
with the dimension of the state space as $O(n^3)$, faster methods that reach an
accurate solution in a small number of iterations are strongly preferred.
We also note that, if one uses the $\ell_1$ regularization penalty solely as a
heuristic for encouraging sparsity, then finding an exact (locally) optimal
solution may be less important than merely finding a solution with a reasonable
sparsity pattern, which can indeed be accomplished by a variety of algorithms.
However, given that we are using the $\ell_1$ heuristic in the first place, and
since in practice the $\ell_1$ penalty has a similar "shrinkage" effect as
increasing the respective $R$ penalty on the controls, it is reasonable to seek
out as accurate a solution as possible to this optimization problem. We
demonstrate that for all levels of
accuracy and on both sets of examples considered, our second order method is
significantly faster than existing approaches. In particular, for
large problems with thousands of states, our method reaches a
reasonable level of accuracy in minutes whereas previous approaches
take hours.

Specifically, we compare our algorithm to two other approaches: the
original alternating direction method of multipliers (ADMM) method
from \cite{lin2013design} which has as an inner loop the
Anderson-Moore method; and iterative soft-thresholding (ISTA), a
proximal gradient approach which
iterates between a gradient step and the soft-thresholding
projection. In the ISTA implementation, in order to ensure
that we maintain the stability of $A + BK$ we perform a line search to
choose the step size for each iteration. Although it is also possible
to add acceleration to ISTA resulting in the FISTA algorithm (which in
theory has superlinear convergence), in practice this performed worse
than ISTA, most likely due to the nonconvexity of the optimization
problem.

In each set of experiments, we first solve the $\ell_1$ regularized
objective with a weight $\lambda$ placed on all elements of $K$. Once
this has converged, we perform a second polishing pass with the
sparsity of $K$ fixed to the nonzero elements of the optimal solution
to the $\ell_1$ problem, optimizing performance on the LQR objective
for a given level of sparsity. The polishing step can also be
performed efficiently using the Newton-CD method and $\Lambda$ with
elements equal to $0$ or $\infty$. Finally,  when solving the $\ell_1$
regularized problem in the first step, we soft-threshold the LQR
solution as described in Section \ref{sec-algorithm}; this is
relatively quick compared to the overall running time and we use the
same initial controller $K^{(0)}$ as the starting point for all
algorithms.

\subsection{Mass-spring system}

In our first example we consider the mass-spring system from
\cite{lin2013design}
describing the displacement of $N$ masses connected on a line. The
state space is comprised of the position and velocity of each mass
with dynamics given by the linear system
\begin{align}
 A = \left[ \begin{array}{rr} 0 & I \\ T & 0 \end{array} \right], &&
 B = \left[ \begin{array}{r} 0 \\ I \end{array} \right]
\end{align}
where $I$ is the $N \times N$ identity matrix, and $T$ is an $N \times N$
tridiagonal symmetric Toeplitz matrix of the form
\begin{equation}
 T = \left[ \begin{array}{rrrr}
-2 &  1 &  0 &  0 \\
 1 & -2 &  1 &  0 \\
 0 &  1 & -2 &  1 \\
 0 &  0 &  1 & -2
\end{array} \right];
\end{equation}
we take $Q = I, R = 10I$, $W = BB^T$ as in the previous paper.

\begin{figure}
\includegraphics{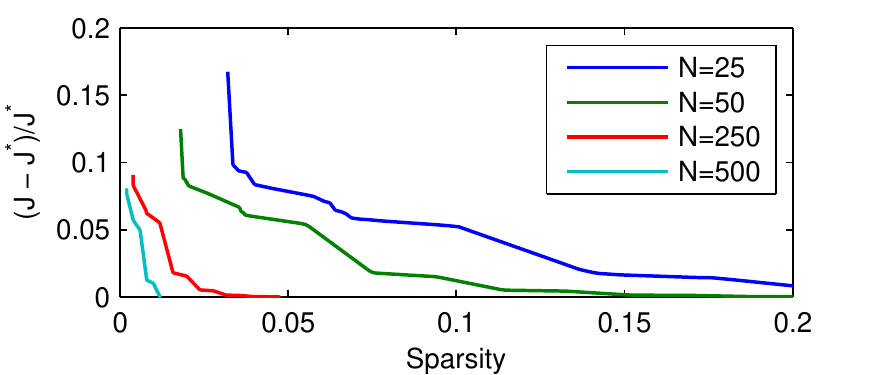}
\caption{Comparison of sparse controllers to the optimal LQR control
  law $J^*$ for varying levels of sparsity on the
  mass-spring system.}
\label{fig-springs-sweep}
\end{figure}

We begin by characterizing the trade-off between sparsity and system
performance by sweeping across 100 logarithmically spaced values of
$\lambda$. For the system with $N=50$ springs, the results shown in
Figure \ref{fig-springs-sweep} are nearly identical to those reported
by \cite{lin2013design}, although their methodology includes an
additional loop and iteratively solving a series of reweighted $\ell_1$
problems. For all systems, the leftmost point represents a control law
based almost entirely on local information---although the results
shown penalize the elements of $K$ uniformly, we also found that by
regularizing just the elements of $K$ corresponding to nonlocal
feedback we were able to find stable local control laws in all
examples. As $\lambda$ decreases, the algorithm finds controllers that
quickly approach the performance of LQR and in the smallest example we
require a controller with 18\% nonzero elements to be within $0.1\%$
of the LQR performance; in the largest example we require only and
4.0\% sparsity to reach this level. This demonstrates the trend that
we anticipate: larger systems require comparatively sparser
controllers for optimal performance.



\begin{figure}
\includegraphics[width=1.72in]{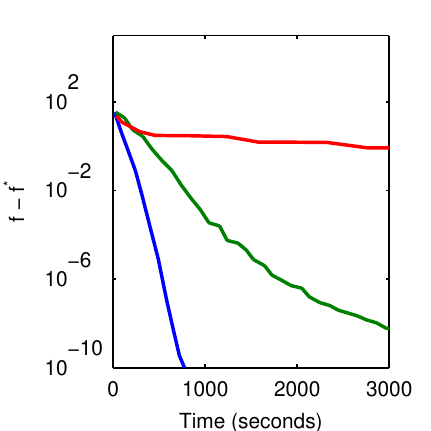}
\includegraphics[width=1.72in]{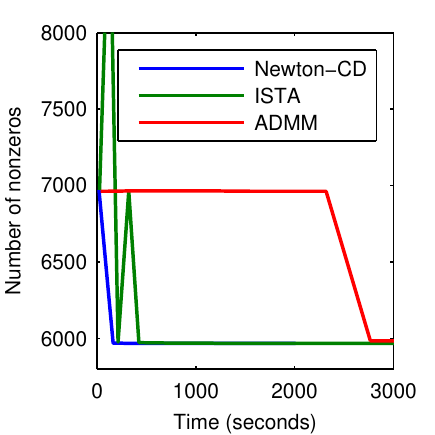}
\includegraphics{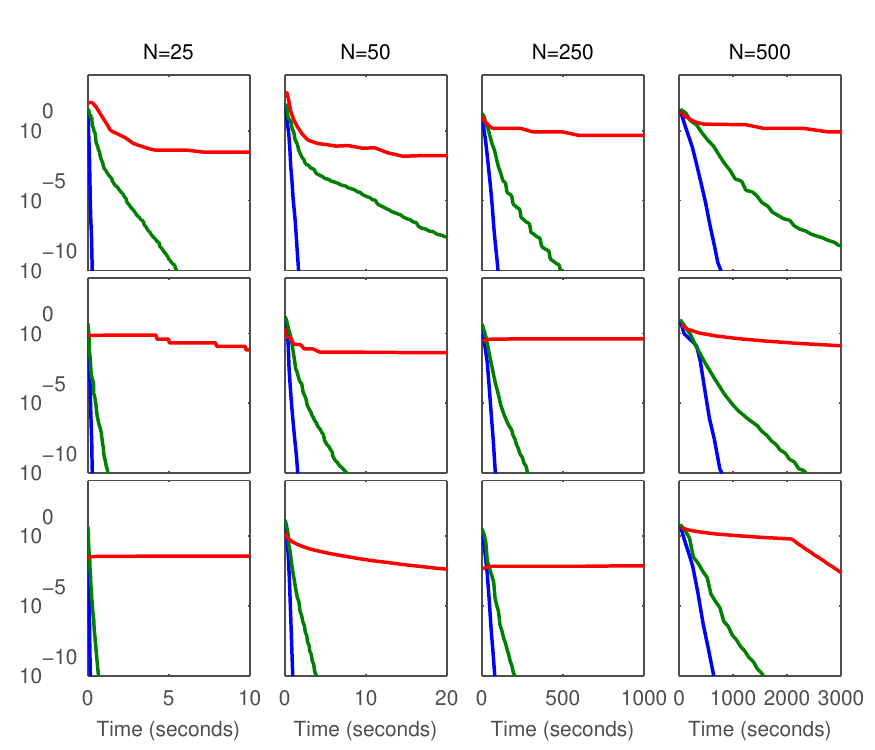}\
\caption{Convergence of algorithms on mass-spring system with $N=500$
  and $\lambda=10$ (top left); the sparsity found by each algorithm
  for the same system (top right); and across many settings
  with one column per example and rows corresponding to different settings of $\lambda$ with
  $\lambda_1 = [10, 10, 1, 1]$, $\lambda_2 = [1, 1, 0.1, 0.1]$ and $\lambda_3 =
  [0.1, 0.1, 0.01, 0.01]$ (bottom).}
\label{fig-springs}
\end{figure}

Next we compare running times of each algorithm by fixing $\lambda$
and considering the convergence of the objective value $f$ at each
iteration to the (local) optimum  $f^*$. Figure \ref{fig-springs}
shows three such fixed $\lambda$
settings corresponding to the levels of sparsity of interest in the
mass-spring system and we in all settings that the Newton-CD method
converges far more quickly than other methods. In the largest
system considered ($N=500$) with $\lambda = 0.1$ (top left), it
converges to a solution accurate to $10^{-8}$ in less than 11 minutes
whereas ADMM has not reached an accuracy of $10^{-1}$ after over two
hours. In addition, the sparsity pattern in the intermediate solutions
do not typically correspond to that of the $\ell_1$ solution, as can
be seen in Figure \ref{fig-springs} (top right). For smaller
examples, ISTA is competitive but as the size of
the system grows, the many iterations that it requires to converge
become more expensive, a behavior that is highlighted in the convergence
for the $N=500$ system (rightmost column). Finally, we note that
Newton-CD performs especially well on $\lambda$ corresponding to
sparse solutions (top row) due to the active set method exploiting
sparsity in the solution.

\begin{figure}
\includegraphics{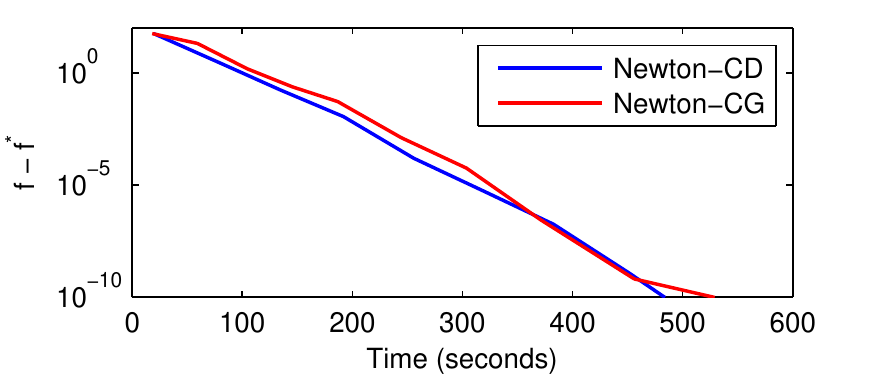}
\caption{Convergence of Newton methods on the polishing step for
  the mass-spring system with $N=500$ and $\lambda=100$.}
\label{fig-polish}
\end{figure}

In addition to solving the $\ell_1$ problem, the Newton-CD method can
also be used for the polishing step of finding the optimal controller
with a fixed sparsity structure. In Figure \ref{fig-polish}, we
compare Newton-CD to the conjugate gradient approach of
\cite{lin2013design} which can be seen as a Newton-Lasso method for
the special case of $\Lambda$ with entries $0$ or $\infty$. Here we
see that performance on the polishing step is comparable with both
methods converging quickly and using the same number of outer loop
iterations. We note that the conjugate gradient approach could also be
extended to work for general $\Lambda$ by using an orthant-based
approach (for example, see \cite{olsen2012newton}), but we do not
pursue that direction in this work.

\subsection{Wide-area control in power systems}

Following \cite{dorfler2013sparsity}, which applied the sparse optimal
control framework and the ADMM algorithm to this same problem, our
next examples consider the task of controlling inter-area oscillations
in a power network via wide-area control.  These examples highlight
the computational benefits of our algorithmic approach even more so than
the synthetic examples above.

To briefly introduce the domain (\cite{dorfler2013sparsity} explains
the overall setup in more detail), we are concerned here with the
problem of frequency regulation in a AC transmission grid.  We employ
a linearized approximation where for each generation the system state
consists of the power angle $\theta_i$, the mismatch between the
rotational velocity and the reference rotational velocity $\omega_i -
\omega^{\mathrm{ref}}$, as well as a number of additional states $x_i$
characterizing the exciters, governors, and/or power system stabilizer
(PSS) control loops at  each generator (typically operating on much
faster time scale).  The system dynamics can be written generically as
\begin{equation}
\begin{split}
\dot{\theta} & = \omega - \omega^{\mathrm{ref}} \\
\dot{\omega} & = (Y_{GG} - Y_{GL} Y_{LL}^{-1} Y_{LG})\theta + f(x)
\end{split}
\end{equation}
where $Y$ is an approximate DC power flow susceptance matrix; $G$ and
$L$ are the generator and load nodes; and $f(x)$
denotes the local control dynamics. Importantly, there is a coupling
between the generators induced by the network dynamics, which can
create oscillatory modes that cannot easily be stabilized by local
control alone.  The control actions available to the system
effectively involve setting the operating points for the inner loops
of the power system stabilizers.

The examples we use here are all drawn from the Power Systems Toolbox,
in particular the MathNetEig package, which provides a set of routines
for describing power networks, generators, exciters and power system
stabilizers, potentially at each generator node, and also has routines
for analytically deriving the resulting linearized systems.  We
evaluate our approach on all the larger examples included with this
toolbox, as well as the New England 39 bus system used in
\cite{dorfler2013sparsity} (which is similar to the 39 bus system
included in the power system toolbox, but which includes power system
stabilizers at 9 of the 10 generators, limits the type of external
control applied to each PSS, and which allows for no control at one of
the generators).  To create a somewhat larger system
than any of those included in the toolbox, we also modify the PST 50
machine system to include power system stabilizers at each node,
resulting in a $n=500$ state system to regulate with our sparse
control algorithm.

\begin{figure}
\includegraphics{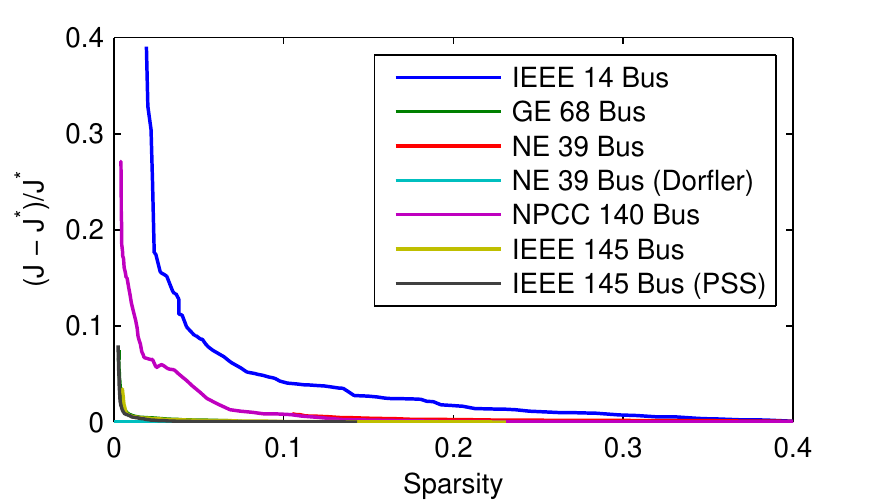}
\caption{Comparison of performance to LQR solution for varying levels of
  sparsity on wide-area control in power networks.}
\label{fig-wac-sweep}
\end{figure}

\begin{figure}
\includegraphics[width=3.5in]{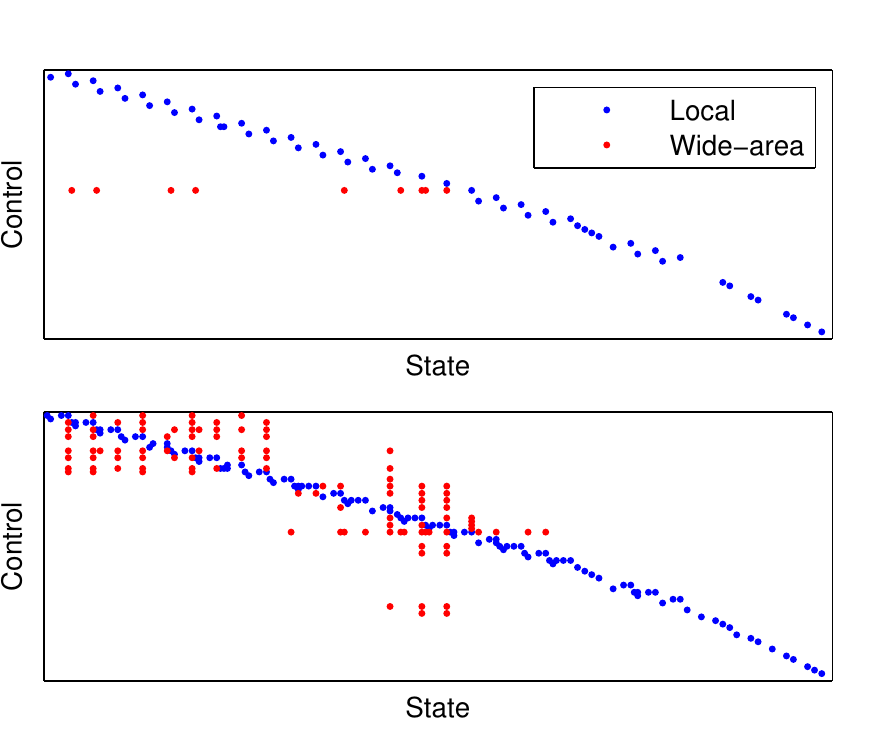}
\caption{Sparsity patterns for wide-area control in the NPCC 140 Bus power
  system: the sparsest stable solution found (top) and the sparsest
  solution achieving performance within 10\% of optimal (bottom)}
\label{fig-wac-comm}.
\end{figure}

As in the previous example, we begin by considering the sparsity/performance
trade-off by varying the regularization parameter $\lambda$, shown in Figure
\ref{fig-wac-sweep}. Here we see that for several powers systems under
consideration, near optimal performance is achieved by an extremely sparse
controller depending almost exclusively on local information. As in
the mass-spring system, in addition to $\Lambda$ with uniform weights,
we also used a structured $\Lambda$ to find local controllers; here
were able to find stable local control laws for every example with the
exception of PST 48. For this power system, we show the sparsity
pattern of the sparsest stable controller in Figure \ref{fig-wac-comm}
along with that of the controller achieving performance within 10\% of
the full LQR optimum. Finally, we note that in general the larger
power systems admit controllers with relatively more sparsity as was
the case in the mass-spring system.

\begin{figure}
\includegraphics[width=1.72in]{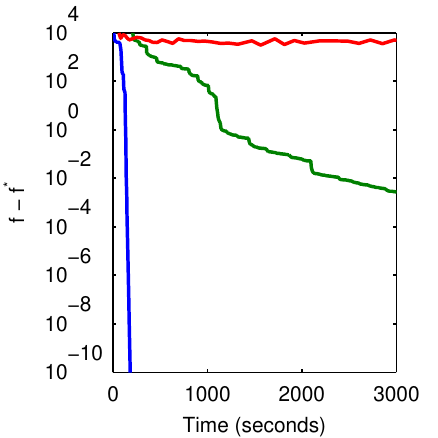}
\includegraphics[width=1.72in]{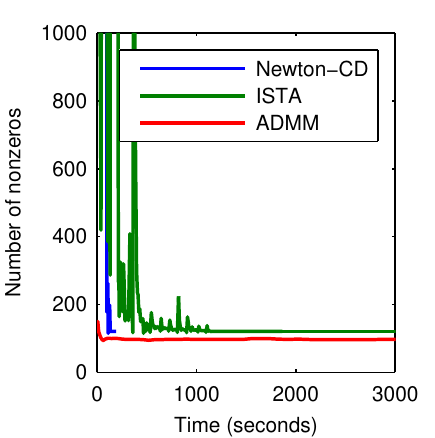}
\caption{Convergence of algorithms on IEEE 145 Bus (PSS) wide-area
  control example with $\lambda = 100$ (left) and the number of
  nonzeros in the intermediate solutions (right).}
\label{fig-ms-large}
\end{figure}

\begin{figure*}
\vspace{-0.9in}
\includegraphics{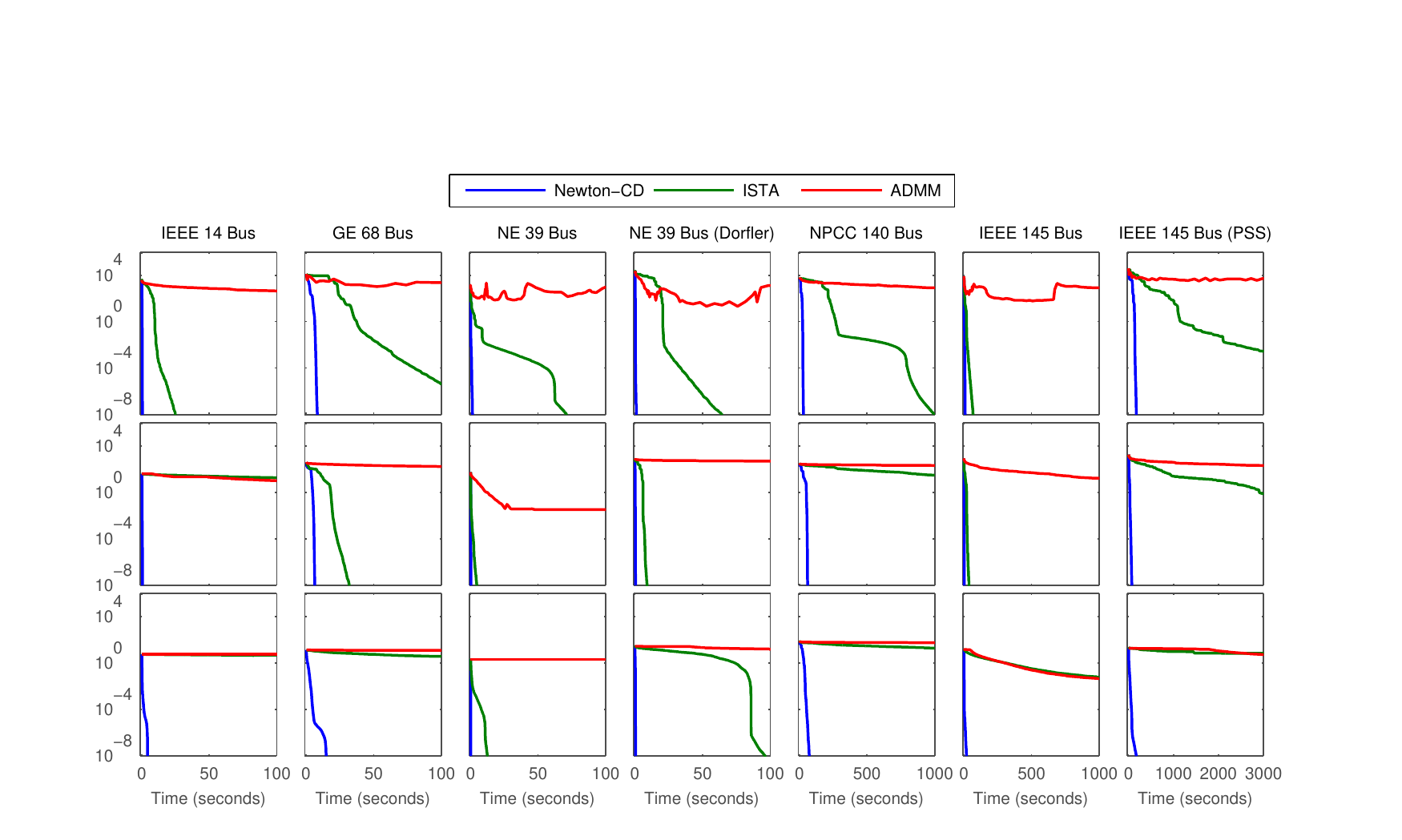}
\caption{Convergence of algorithms on wide-area control across all
  power systems with three choices of $\lambda$ corresponding to
  performance within $10\%$, $1\%$ and $0.1\%$ of LQR. Columns
  correspond to power systems and rows correspond to different
  choices of $\lambda$ with largest on top.}
\label{fig-ms}
\end{figure*}

Computationally, we consider the convergence on the largest power
system example in Figure \ref{fig-ms-large} and observe a dramatic
difference between Newton-CD and previous algorithms: Newton-CD has
converged to an accuracy better than $10^{-8}$ in less than 173
seconds while ADMM is not within $10^3$ after over an hour. In
addition, the sparsity pattern of the intermediate solutions found by
ADMM is significantly different than that of Newton-CD and ISTA which
have converged to the $\ell_1$ regularized solution with much higher
accuracy. In Figure
\ref{fig-ms} we consider convergence across all power systems with
three values of $\lambda$ chosen such that the resulting controllers
have performance within $10\%$, $1\%$ and $0.1\%$ of
LQR. Here we see similar results as in the large system with Newton-CD
converging faster across all examples and choices of $\lambda$ and the
differences being orders of magnitude in many cases. We note that
algorithm benefits significantly from a high level of sparsity for
these choices of $\lambda$ since for most power systems considered,
the controller achieving performance within $0.1\%$ of LQR is still
quite sparse.

\section{Conclusion}
\label{sec-conclusion}

In this paper we develop a fast second order algorithm for the sparse
LQR problem with the goal of designing sparse control laws for  large
distributed systems with thousands of nodes. Intuitively we expect
that distributed systems characterized by sparsity in the dynamics
(e.g. neighbor interactions in the mass-spring system and the graph Laplacian in
power networks) can be well controlled with only limited information
sharing; we see in the experimental results that as these
systems increase in size, they are controllable by laws that
become increasingly sparse. Computationally, the design of
efficient algorithms is complicated by both the Lyapunov equations which
give rise to $O(n^3)$ complexity, and the nonsmooth $\ell_1$
penalty. Our work limits this complexity through an efficient
Newton-Lasso algorithm reducing the time required to solve the
$\ell_1$ problem to that of previous results with fixed
structure. However, the $O(n^3)$ complexity inherent in the Lyapunov
equations pose a significant bottleneck to scaling to systems beyond
thousands of nodes and thus decomposition methods allowing smaller
subproblems to be solved independently are an interesting direction
for future research.

\bibliographystyle{IEEEtran}
\bibliography{sparseoptcontrol}

\end{document}